\documentclass[11pt,a4paper]{article}%
\usepackage[centertags]{amsmath}
\usepackage{amsfonts}
\usepackage{amssymb}
\usepackage{amsthm}
\usepackage{epsfig}
\usepackage{setspace}
\usepackage{ae}
\usepackage{eucal}
\usepackage[usenames]{color}%
\usepackage{graphicx}

\theoremstyle{plain}
\newtheorem{thm}{Theorem}[section]

\newtheorem{cor}[thm]{Corollary}
\newtheorem{prop}[thm]{Proposition}
\theoremstyle{definition}

\theoremstyle{remark}

\newtheorem{rem}[thm]{Remark}

\begin{document}

\title{Global analytic expansion of solution for a class of linear parabolic systems with coupling of first order derivative terms  }
\author{J\"org Kampen }
\maketitle

\begin{abstract}
We derive global analytic representations of fundamental solutions for a class of linear parabolic systems with full coupling of first order derivative terms where coefficients may depend on space and time. Pointwise convergence of the global analytic expansion is proved. This leads to analytic representations of solutions of initial-boundary problems of first and second type in terms of convolution integrals or convolution integrals and linear integral equations. The results have both analytical and numerical impact. Analytically, our representations of fundamental solutions of coupled parabolic systems may be used to define generalized stochastic processes.
Moreover, some classical analytical results based on a priori estimates of elliptic equations are a simple corollary of our main result. Numerically, accurate, stable and efficient schemes for computation and error estimates in strong norms can be obtained for a considerable class of Cauchy- and initial-boundary problems of parabolic type. Important instances of application are representations of solutions of multidimensional Burgers equations with forcing and potential initial conditions and Pauli equation describing the non-relativistic limit of Dirac theory for electrons in a magnetic field.

Warning: The current analysis holds only in special cases, essentially in the case of scalar equations. A more involved analytic expansion is necessary (and possible) for systems and will be communicated soon.
\end{abstract}

\footnotetext[1]{Weierstrass Institute for Applied Analysis and Stochastics,
Mohrenstr. 39, 10117 Berlin, Germany.
\texttt{{kampen@wias-berlin.de}}.}

\section{Introduction}
We consider linear equations of the form

\begin{equation}\label{parasyst1}
\frac{\partial \mathbf{u}}{\partial t}= \nabla^2  \mathbf{u}
+\mathbf{B}\cdot\mathbf{\nabla u}
\end{equation}
on a domain $D=(0,T]\times \Omega$ with $\Omega \subseteq {\mathbb R}^n$ a bounded domain, and where
\begin{equation}
\mathbf{u}=(u_1,\cdots ,u_n)^T
\end{equation}
is a vector-valued function and 
$\mathbf{B}=(B^1,\cdots,B^n)$ is an $n$-tuple of matrix-valued functions $B^i=(b^i_{jk})$ where each entry $b^i_{jk}$ possibly depends on space and time.
More precisely, we understand $\mathbf{B}\cdot\mathbf{\nabla u}$ as a vector the $i$th component of which is given by
\begin{equation}
(\mathbf{B}\cdot\mathbf{\nabla u})_i:=\sum_{j,k=1}^n b^i_{jk}\frac{\partial u_j}{\partial x_k},
\end{equation}
such that general linear coupling of first order terms can be expressed.
This means that in coordinates equation (\ref{parasyst1}) is given by
\begin{equation}\label{parasyst2}
\frac{\partial u_i}{\partial t}=\sum_{j=1}^n \frac{\partial^2 u_i}{\partial x_j^2} 
+\sum_{j,k=1}^n b^i_{jk}\frac{\partial u_j}{\partial x_k}
\end{equation}
for $1\leq i\leq n$. We are looking for an analytic representation of the solution $(t,x,s,y)\rightarrow {\bf p}(t,x;s,y)$ for (\ref{parasyst1}), (\ref{parasyst2}) with Dirac distributions $\delta_y(x)=\delta(x-y)$ as initial conditions, i.e. for a representation of the fundamental solution. In the time-homogenous case dependence of time is only dependence of $t-s$, so that in this case we also write the fundamental solution in the form $(t,x,y)\rightarrow {\bf p}(t,x;y)$. For our representations of the fundamental solution for equations with time-dependent coefficients we also fix the parameter $s=0$ and write the fundamental solution in the form  $(t,x,0,y)\rightarrow {\bf p}(t,x;0,y)$ for simplicity of notation.  
\begin{rem}
In the following we denote the fundamental solution of a system by bold face letters and use usual type of letters for the fundamental solution of a scalar equation.
\end{rem}

\begin{rem}
In general for parabolic systems of order $2p$ of form
\begin{equation}\label{gensyst}
\frac{\partial u_i}{\partial t}=\sum_{j=1}^N \sum_{|\alpha|\leq 2p}A^{ij}_{\alpha}(t,x)\partial^{\alpha}_x u_j
\end{equation}
(with some natural number $N$) the fundamental solution (or fundamental matrix) $(t,x,s,y)\rightarrow {\bf p}(t,x,s,y)$ is a $N\times N$-matrix of functions on $(0,T]\times\Omega\times (0,T]\times\Omega$ which satisfies
(\ref{gensyst}), and such that 
\begin{equation}\label{cond6}
\lim_{t\downarrow s}\int_{\Omega}f(y){\bf p}(t,x,s,y)dy=f(x)
\end{equation}
for all continuous functions $y\rightarrow f(y)$ in $\overline{\Omega}$. The latter condition is equivalent to the rule that $p(0,x;0,y)=\delta_y(x)=\delta(x-y)$. Here in the general case with coupling of the higher order derivatives $2p$, $p\geq 1$, a vectorial representation of the fundamental solution is not known. It is one of the observations of this paper that a vectorial representation is possible in the case $p=1$, if the only coupling occurring is that via first order terms. Note that $N\neq n$ is possible. Our restriction to the case $N=n$ is not essential but only related to an economy in the use of symbols.
\end{rem}
\begin{rem}
Probabilistic representations of the solution of the Cauchy problem for linear parabolic systems have been obtained first by Stroock (as far as I know, cf. \cite{S}). Note that in \cite{S} a representation is obtained via a fundamental matrix representation different to our representation which is vectorial. Moreover, our approach may be used in order to extend the Feynman-Kac formalism slightly beyond the class of processes considered in \cite{S}, because we may have $n$ diffusion matrices
$a^i_{jk},~~1\leq i\leq n$ in our generalization mentioned in the last section below.
\end{rem}

It turns out that results in the case of time-homogenous coefficients can be extended to the case of time-inhomogeneous coefficients but it is worth to consider the time-homogeneous case separately, because less assumptions have to be made. For this reason we shall state our main theorem in the time-homogenous case separately, i.e. where coefficient functions are of the form $x\rightarrow b^i_{jk}(x)$.
We shall assume that the functions $x\rightarrow b^i_{jk}(x)$ and their derivatives are uniformly bounded by powers of a generic constant $c$ such that
\begin{equation}\label{coeff}
|\partial_x^{\alpha}b^i_{jk}|\leq c^{|\alpha|}
\end{equation}
for all multiindices $\alpha=(\alpha_1,\cdots ,\alpha_n)$. Here $\partial^{\alpha}_x=\frac{\partial^{|\alpha|}}{\partial x_1^{\alpha_1}\cdots\partial x_n^{\alpha_n}}$ denotes the partial derivative operator of order $\alpha$ with respect to $x=(x_1,\cdots,x_n)$.
If the coefficients are time-dependent functions of form $(t,x)\rightarrow b^i_{jk}(t,x)$, then we shall assume in addition that
\begin{equation}\label{coeft}
|\partial_t^{m}b^i_{jk}|\leq C m!~~\mbox{for all integers $m\geq 0$,}
\end{equation}
all $1\leq i,j,k\leq n$, $\partial_t^{m}=\frac{\partial^m}{\partial t^m}$ is the derivative of order $m$ with respect to time.
 
Note that assumption \ref{coeff} holds for arbitrary finite Fourier series. Clearly it holds also on a bounded domain for polynomial coefficients, and on such domains multivariate polynomials can approximate all continuous functions in the supremum norm. This means that we are flexible enough for numerical applications. Indeed the treatment for higher dimensional scalar parabolic problems based on analytic expansions of the type considered here (simplified to the scalar case) showed accurate and fast computations as well as error estimates in strong norms (cf. \cite{KKS}). The main reasons, however, to introduce assumption \ref{coeff} are analytical. First assumption \ref{coeff} implies that $b^i_{jk}$ are globally analytic, i.e. for all $y\in {\mathbb R}^n$ $b^i_{jk}$ equals its Taylor expansion, i.e. we have
\begin{equation}\label{powerb}
b^i_{jk}(x)=\sum_{|\alpha|\geq 0} \frac{b^i_{jk\alpha}(y)}{\alpha!}\Delta x^{\alpha},
\end{equation}
where $\alpha$ is a multiindex and $\Delta x=(x-y)$, and $b^i_{jk\alpha}(y)=\partial^{\alpha}_x b^i_{jk}(y)$. This makes it possible to write down explicit solutions of parabolic systems of type \ref{parasyst1} in terms of power series of type \ref{powerb}. Second, the proof indicates that the assumptions made here cannot be weakened in general.

Essential parts of our considerations can be generalized to parabolic systems with space-time dependent second-order terms, i.e. equations of the form
\begin{equation}\label{parasyst3}
\frac{\partial u_i}{\partial t}=\sum_{j,k=1}^n a^i_{jk}\frac{\partial^2 u_i}{\partial x_j\partial x_k} 
+\sum_{j,k=1}^n b^i_{jk}\frac{\partial u_j}{\partial x_k},
\end{equation}
where the scalar functions $a^i_{jk}$ may depend on space and time. It turns out that the convergent analytical solutions (in case (\ref{parasyst2})) are building blocks for the representations of solutions of parabolic systems with potential and source terms.
As examples, let us consider a Cauchy problem and a standard initial-boundary problem which occur in the vector-valued as well as in the scalar case. The finite horizon Cauchy problem for parabolic systems of type (\ref{parasyst2}) is
\begin{equation}\label{cp}
\left\lbrace \begin{array}{ll}
\frac{\partial u_i}{\partial t}-\sum_{j=1}^n \frac{\partial^2 u_i}{\partial x_j^2} 
-\sum_{jk} b^i_{jk}\frac{\partial u_j}{\partial x_k}=f_i~~\mbox{in ${\mathbb R}^n\times (0,T]$}\\
\\
u_i(0,x)=\phi_i(x)~~\mbox{on~${\mathbb R}^n$},
\end{array}\right.,
\end{equation}
where $T>0$ and for $1\leq i \leq n$.
\begin{rem}
The proof of the pointwise valid representation of the fundamental solution given is valid for bounded domains $\Omega$ and cannot be directly generalized to unbounded domains. This is no essential restriction for numerical treatment, however. Analytically, a generalization is possible, if one considers a slightly different representation (cf. section 4.2.). However, the recursive relations of the expansion coefficients are more complicated and the convergence proof is more involved. Therefore we restrict ourselves to the case of bounded domains $\Omega$ in this paper.
\end{rem}

Another example is the initial-boundary problem of second type. We consider it in the scalar case here. Consider a domain $\Omega\subset {\mathbb R}^n$ and denote the three constituents of boundary of the cylinder by $\Omega_0:=\left\lbrace (t,x)|t=0~\&~x\in\Omega\right\rbrace$  $\Omega\times (0,T)$ by $\Omega_T=\left\lbrace (t,x)|t=T~\&~x\in\Omega \right\rbrace$ and $B=\left\lbrace (t,x)|t\in (0,T)~\&~x\in \partial\Omega \right\rbrace$, where $\partial \Omega$ denotes the boundary of $\Omega$. The initial-boundary problem is of the form 
\begin{equation}\label{ibst}
\left\lbrace \begin{array}{ll}
\frac{\partial u}{\partial t}-\sum_{j=1}^n \frac{\partial^2 u}{\partial x_j^2} 
-\sum_{k} b_{k}\frac{\partial u}{\partial x_k}=f~~&\mbox{in $\Omega\times (0,T]$}\\
\\
u(0,.)=\phi(.)~~&\mbox{on~$\Omega$}\\
\\
\frac{\partial u}{\partial t}+\alpha u=\psi~~&\mbox{on~$B$},
\end{array}\right.
\end{equation}
where $\alpha, \phi,$ and $\psi$ may depend on space and time. With an explicit representation of the fundamental solution we can represent the solution of (\ref{cp}) in terms of convolutions of the initial data and the source data with the fundamental solution, and the solution of (\ref{ibst}) in terms of convolutions of initial data, source data, and a function which is solution of a linear integral equation. It is clear that such representations lead to accurate schemes which have obvious advantages compared to finite difference schemes and other standard schemes. 
\begin{rem}
In (\ref{parasyst2}) we may add potential terms of form $c_i u$ with a coefficient functions $c_i$ which may depend on space and time. Theorem 1 below can be trivially extended to this case. Hence, in 
equations (\ref{cp}) and (\ref{ibst}) we may also add potential terms of form $c_i u$ and representations of solutions in terms of convolutions and linear integral equations (in case of the initial-boundary problem (\ref{ibst}) can be obtained.
\end{rem}
\begin{rem}
More general cylinder domains $D=\cup_{0\leq t\leq T} \Omega_t$ with varying $\Omega_t$ may be considered, of course. 
\end{rem}
This is the first paper on globally pointwise valid analytic expansions of parabolic systems. In the case of scalar equations there are some investigations and applications to problems of computation recently (\cite{KKS} and references). Our result has direct applications to case of the scalar equations, of course. Further comments on the relation to results in the scalar case can be found in Section 7.

The outline of this paper is as follows. In the next section we state the main results concerning the representation of the fundamental solution. In Section 3 we formally compute the analytic expansion of the solution and in Section 4 we prove the pointwise convergence of the analytic representation in the time-homogenous case for a certain limited time horizon $0\leq t\leq T_0$. In Section 5 we extend the results of the preceding Sections to the case where the coefficients may depend on space and time and we show the global convergence for any time horizon $0<T<\infty$. In Section 6 we consider the implications for representations of solutions Cauchy problems and second initial-value boundary problem and briefly discuss the advantages for building efficient numerical schemes. In Section 7 we state some generalizations with general (but uncoupled) diffusion coefficients and show that a result by Varadhan is a rather immediate consequence of our main theorem. We also discuss possible other applications (for example the definition of generalized processes) and give some further comments and an outlook.
\section{Main results on linear parabolic systems}

Since the second order derivative terms in (\ref{parasyst2}) are uncoupled, we may expect that a vectorial representation of the fundamental solution ${\bf p}=(p_1,\cdots ,p_n)$ (instead of an $n\times n$ fundamental matrix) is possible. The natural candidate for such a representation (in the time-homogeneous case) is
\begin{equation}\label{pordj}
p_j(t,x,y)=\frac{1}{\sqrt{4\pi t}^n}\exp\left(-\frac{\sum_{i=1}^n \Delta x_i^2}{4t}+\sum_{k=0}^{\infty} c^j_{k}(x,y)t^k \right),  
\end{equation}
for $j=1,\cdots ,n$, and  for $(t,x) \in(0,T)$, where $\Omega \subseteq {\mathbb R}^n$. 
Here the $c^j_k$ are coupled coefficient functions which are defined explicitly via recursion. For each $j$ the coefficients $c^j_k$ will be defined recursively in terms of function $c^r_{l}$ and their derivatives, where $0\leq l,r\leq k-1$. They are solutions of first order partial differential equations which can be solved explicitly and can be represented in terms of recursively defined power series under the assumption (\ref{powerb}). We shall show that (\ref{pordj}) is valid on some domain $\Omega\times(0,T_0]$. Since it is desirable to have a representation which holds on the whole a domain $\Omega\times(0,T]$ with arbitrary time $T\in (0,\infty)$, in our main theorem we shall consider global representations of an equivalent problem, where the equivalence is via the time transformation
$\tau(t):[0,\infty)\rightarrow [0,1)$ with
\begin{equation}\label{timetrans}
\tau =(1-e^{-\frac{t}{\beta}}),~~\mbox{or}~~t=t(\tau)=-\beta\ln(1-\tau).
\end{equation}
This introduces a time-dependence in the related coefficients 
$c^j_{k,\beta,\tau}$, even in the case of time-homogeneous coefficient functions  $x\rightarrow b^i_{jk}(x)$ in (\ref{parasyst2}).
The main result for parabolic systems of type (\ref{parasyst1}) is formulated in the time-homogenous case, i.e. when the coefficients $b^j_{lm}$ depend only on the spatial coordinates. The extension to the time-dependent case is then the content of the subsequent corollary. 
\begin{thm}
Given assumption (\ref{coeff}) and some domain $\Omega\times (0,T]$ for any finite $T>0$ and any domain $\Omega\subseteq {\mathbb R}^n$ there exist $\beta,\tau>0$ such that the fundamental solution of 
\begin{equation}\label{parasystthm}
\frac{\partial u_i}{\partial \tau}=\frac{\beta}{1-\tau}\sum_{j=1}^n \frac{\partial^2 u_i}{\partial x_j^2} 
+\frac{\beta}{1-\tau}\sum_{j,k=1}^n b^i_{jk}\frac{\partial u_j}{\partial x_k}
\end{equation}
equivalent to (\ref{parasyst1}) (or (\ref{parasyst2})) via (\ref{timetrans})
has the pointwise valid representation
\begin{equation}\label{pj}
p^{\beta ,\tau}_j(\tau,x,0,y)=\frac{1}{\sqrt{4\pi t(\tau)}^n}\exp\left(-\frac{\sum_{i=1}^n \Delta x_i^2}{4t(\tau)}\right)\exp\left(\sum_{k=0}^{\infty} c^j_{k,\beta,\tau}(\tau,x,y)\tau^k \right),  
\end{equation}
for $j=1,\cdots ,n$, and for for $(t(\tau),x) \in(0,T)\times \Omega$, i.e. $\tau \in (0,1-e^{-\frac{T}{\beta}})$, where $\Omega \subseteq {\mathbb R}^n$. 
For the coefficient functions $c^j_k$ the following holds: for $k=0$ we have
\begin{equation}\label{c0}
 c^j_{0,\beta,\tau}(\tau,x,y)= c^j_0(x,y)=\sum_m (y_m-x_m)\int_0^1 \sum_lb^j_{l,m}(y+s(x-y))ds,
\end{equation}
and for all $k\geq 1$ we have
\begin{equation}
c^j_{k,\beta,\tau}(\tau, x,y)=\int_0^1 R^j_{k-1,\beta,\tau}(t,y+s(x-y),y)s^{k\frac{1-\tau}{\beta}-1}ds
\end{equation}
with 
\begin{equation}\label{tk}
\begin{array}{ll}
R^j_{k-1,\beta,\tau}(t,x,y)=&\frac{\partial}{\partial \tau}c^j_{k-1,\beta,\tau}+\Delta c^j_{k-1,\beta,\tau}+\sum_{l=1}^n\sum_{r=0}^{k-1}\left( \frac{\partial}{\partial x_l}c^j_{r,\beta,\tau}\frac{\partial}{\partial x_l}c^j_{k-1-r,\beta,\tau}\right)\\
\\
&+\sum_{lm} b^j_{lm}(x)\frac{\partial}{\partial x_m}c^l_{k-1}
\end{array}
\end{equation}
More explicitly, we have
\begin{equation}\label{powerco}
\begin{array}{ll}
 c^j_{0,\beta,\tau}(\tau,x,y)=c^j_0(x,y)=&-\sum_{l,m} \sum_{\gamma}b^{j}_{lm\gamma}(y)\Delta x^{\gamma+1_i}\frac{1}{1+|\gamma|}\\
\\
&\equiv \sum_{\gamma} c^j_{0\gamma}\Delta x^{\gamma}
\end{array}
\end{equation}
and, given the power series representation 
\begin{equation}\label{powbeta}
c^j_{k-1,\beta,\tau}(\tau,x,y)
=\sum_{\gamma,l}c^j_{(k-1)\gamma l}(y)\Delta x^{\gamma}\tau^l
\end{equation}
we have
\begin{equation}\label{powerck}
\begin{array}{ll}
 c^j_{k,\beta,\tau} (\tau, x,y) =\sum_{\gamma,l}lc^j_{(k-1)\gamma l}(y)\Delta x^{\gamma}t^l+\\
\\
 \sum_{\gamma}{\big \{}\sum_i\sum_{\rho + \alpha = \gamma} (\rho_i +1)(\alpha_i +1 )c^j_{r(\beta +1_i)}c^j_{(k-1-r)(\alpha +1_i)}   \\
\\
 +\sum_i (\gamma_i+2)(\gamma_i+1)c_{k(\gamma+2_i)}+ 
\sum_{\rho + \alpha = \gamma}(\sum\frac{1}{\beta !}b^{j}_{lm,\rho}(y)\times\\
\\
(\alpha_i+1)c_{(k-1)(\alpha+1_i)}{\big \}}\left( \sum_{\delta=0}^{\gamma}p_{k\delta}^{y\gamma} \Delta x^{\delta}\right),  
\end{array}
\end{equation}
where with $\displaystyle \delta_{\Sigma} := \sum_{i=1}^n \delta_i$, and
\begin{eqnarray}\label{abinom}
\sum_{\delta=0}^{\gamma}p_{k\delta, \beta,\tau}^{y\gamma} \Delta x^{\delta}&=&\sum_{\delta=0}^{\gamma}\frac{\beta}{(1-\tau)\delta_{\Sigma} +k}\nonumber\\
\\
&\times&\left[  \prod_{i=1}^n\left( \frac{\gamma_i !}{\delta_i ! (\gamma_i - \delta_i)!}\right)  y^{(\gamma - \delta)}\right] \Delta x^{\delta}.\nonumber
\end{eqnarray}
\end{thm}
\begin{rem}
In (\ref{powbeta}) the notation
\begin{equation}
c^j_{k,\beta,\tau}(x,y)=\sum_{\gamma l}c^j_{(k-1)\gamma l,\beta,\tau}(y)\Delta x^{\gamma}\tau^l
\end{equation}
may be expected, but we dropped the $\beta,\tau$ indices in order to keep some notational simplicity in (\ref{powerck}).
\end{rem}
\begin{rem}
Note that $c^j_{0,\beta,\tau}$ does not depend on $\tau$ (and $\beta$). This means that in (19) $\frac{\partial}{\partial t}c_{k-1}$ differs from $0$ only for $k\geq 2$. 
\end{rem}

\begin{cor}
Consider the same situation as in the preceding theorem, except that there are time dependent coefficient functions $(t,x)\rightarrow b^i_{jk}(t,x)$ and in addition (\ref{coeft}) holds. Then a analogous statement as in the preceding theorem holds with recursive . 
\end{cor}

\section{Formal computation of solution of parabolic systems of type \eqref{parasyst1} }
First we consider the equation (\ref{parasyst1}) (or, equivalently, (\ref{parasyst2})) without the time transformation (\ref{timetrans}), and with time-homogenous coefficients, i.e. where the coefficient functions $x\rightarrow b^i_{jk}(x)$ depend only on the spatial variable $x$. We consider the ansatz
\begin{equation}\label{pjo}
p_j(t,x,y)=\frac{1}{\sqrt{4\pi t}^n}\exp\left(-\frac{\sum_{i=1}^n \Delta x_i^2}{4t}+\sum_{k=0}^{\infty} c^j_k(x,y)t^k \right).  
\end{equation}
We derive recursive relations for the coefficients $c_k$. In a second step, assuming (\ref{coeff}), and therefore global analyticity of the $b^i_{jk}$, we derive the explicit solution in terms of Taylor power series of $b^i_{jk}$. 
For the time derivative we get
\begin{equation}
\frac{\partial p_j}{\partial t}(t,x)=\left(-\frac{n}{2t}+\frac{\sum_i \Delta x_i^2}{4t^2}+\sum_k kc^j_k(x,y)t^{k-1}\right)p_j (t,x,y).
\end{equation}
For the first and second spatial derivatives we get
\begin{equation}
\frac{\partial p_j}{\partial x_l}=\left(\frac{-\Delta x_l}{2t}+\sum_k \frac{\partial}{\partial x_l}c^j_k(x,y)t^k \right) p_j(t,x,y),
\end{equation}
and
\begin{equation}
\begin{array}{ll}
\frac{\partial^2 p_j}{\partial x_l^2}=&{\Bigg (}-\frac{1}{2t}+\sum_k \frac{\partial^2}{\partial x_l^2}c^j_k(x,y)t^k \\
\\
&+\left(-\frac{\Delta x_l}{2t}+\sum_k \frac{\partial}{\partial x_l}c^j_k(x,y)t^k \right)^2{\Bigg )} p_j(t,x,y).
\end{array}
\end{equation}
Plugging into (\ref{parasyst2}) and ordering with respect to the terms $t^{-2},t^{-1}$ etc. we get the following
recursive relations for the $c^j_k$, where $1\leq j\leq n$:
\begin{equation}\label{tsmminus2}
t^{-2}:~~\frac{\sum_i \Delta x_i^2}{4t^2}=\sum_{l}\frac{ \Delta x_l^2}{4t^2}
\end{equation}

\begin{equation}\label{tsmminus1}
t^{-1}:~~-\frac{n}{2t}=-\sum_{l}\frac{1}{2t}-\frac{1}{2t}\left( \sum_l \Delta x_l\frac{\partial c_0^j}{\partial x_l}-\sum_{lm} b^j_{lm}(x)\Delta x_m\right) , 
\end{equation}
and for all $k-1\geq 0$.
\begin{equation}\label{tsmk}
\begin{array}{ll}
t^{k-1}:~~kc^j_k+\sum_l \Delta x_l\frac{\partial c_k^j}{\partial x_l}=&\Delta c^j_{k-1}+\sum_{l=1}^n\sum_{r=0}^{k-1}\left( \frac{\partial}{\partial x_l}c^j_r\frac{\partial}{\partial x_l}c^j_{k-1-r}\right)\\
\\
&+\sum_{lm} b^j_{lm}(x)\frac{\partial}{\partial x_m}c^l_{k-1}\equiv R^j_{k-1}(x,y).
\end{array}
\end{equation}
Note that the first order coupling of the system is essentially reflected in the recursive first order partial differential equations starting from (\ref{tsmk}). This would be different if we had coupling via the second order terms and it makes the solution of the system much easier.
Note that equation (\ref{tsmminus2}) is satisfied. Equation (\ref{tsmminus1}) is equivalent to
\begin{equation}
\sum_l \Delta x_l\frac{\partial c_0^j}{\partial x_l}=-\sum_{l,m} b^j_{lm}(x)\Delta x_m, 
\end{equation}
with the solution
\begin{equation}
 c^j_0(x,y)=\sum_m (y_m-x_m)\int_0^1 \sum_lb^j_{l,m}(y+s(x-y))ds
\end{equation}
and for all $k\geq 1$ we have
\begin{equation}
c^j_k(x,y)=\int_0^1 R_{k-1}(y+s(x-y),y)s^kds
\end{equation}
with $R_{k-1}$ as in equation \eqref{tk}.
Next we compute the solution explicitly doing the integral for $c_0$ first. We abbreviate $\Delta x=(x-y)$ with components $\Delta x_i=(x-y)_i$ and for a multiindex $\alpha=(\alpha_1,\cdots, \alpha_n)$ we write $\Delta x^{\alpha}:=\Pi_{i=1}^n\Delta x_i^{\alpha_i}$. Furthermore, we define $|\alpha|=\sum_i \alpha_i$ If
\begin{equation}
b^j_{lm}(x)=\sum_{\gamma} \frac{1}{\gamma !}b^{j}_{lm,\gamma}(y)(\Delta x)^{\gamma},
\end{equation}  
along with some multiindex $\gamma$, then
\begin{equation}
\begin{array}{ll}
 c^j_0(x,y)&=-\sum_m \Delta x_m \int_0^1 \sum_l b^j_{lm}(y + s\Delta x)ds\\
\\
&=-\sum_{l,m} \Delta x_m \int_0^1 \sum_{\gamma} b^{j}_{lm\gamma}(y)(s\Delta x)^{\gamma}ds\\
\\
&=-\sum_{l,m} \Delta x_m \sum_{\gamma} b^{j}_{lm\gamma}(y)\Delta x^{\gamma} \int_0^1 s^{|\gamma|}ds\\
\\
&=-\sum_{l,m} \sum_{\gamma}b^{j}_{lm\gamma}(y)\Delta x^{\gamma+1_i}\frac{1}{1+|\gamma|}s^{|\gamma|+1}\Big |_0^1\\
\\
&=-\sum_{l,m} \sum_{\gamma}b^{j}_{lm\gamma}(y)\Delta x^{\gamma+1_i}\frac{1}{1+|\gamma|}\\
\\
&\equiv \sum_{\gamma} c^j_{0\gamma}\Delta x^{\gamma}.
\end{array}
\end{equation}
Next we compute $c^j_k$ for $k\geq 1$.  We have
\begin{equation}
\begin{array}{ll}
 c^j_{k} (x,y) =&\int_0^1 \{\sum_i \sum_{r=0}^{k-1} \dfrac{\partial c^j_r}{\partial x_i}\dfrac{\partial c^j_{k-1-r}}{\partial x_i}\\
\\
 &+ \Delta c_{k-1} + \sum_i b^j_{lm} \dfrac{\partial c_{k-1}}{\partial x_i}\} (y+s(x-y))s^{k-1}ds.
\end{array}
\end{equation}
Assuming that $c^ j_{k-1}$ equals its Taylor series for every $y\in {\mathbb R}^n$, i.e.
\begin{equation}
c^j_{k-1}(x)=\sum_{\gamma} c^j_{(k-1)\gamma}(y)\Delta x^{\gamma},
\end{equation}
then we may evaluate the derivatives occurring in $R_{k-1}$ as follows:
\begin{equation}
 \dfrac{\partial c^j_{k-1}}{\partial x_i}=\sum_{\gamma}(\gamma_i+1)c_{(k-1)(\gamma+1_i)}(y)\Delta x^{\gamma},
\end{equation}
and
\begin{equation}
\dfrac{\partial^2 c^j_{k-1}}{\partial x_i^2}=\sum_{\gamma}(\gamma_i+2)(\gamma_i+1)c_{k(\gamma+2_i)}(y)\Delta x^{\gamma},
\end{equation}
and
\begin{equation}
 \dfrac{\partial c^j_{r}}{\partial x_i}\dfrac{\partial c^j_{k-1-r}}{\partial x_i}=\sum_{\gamma}\left\lbrace \sum_{\beta + \alpha = \gamma} (\beta_i +1)(\alpha_i +1 )c^j_{r(\beta +1_i)}c^j_{(k-1-r)(\alpha +1_i)} \right\rbrace \Delta x^{\gamma}. 
\end{equation}
For the multiindex $\gamma$, we have
\begin{eqnarray}\label{abinom}
 P_k^{\gamma}(x,y) &:=&\int_0^1 (y + s(x-y))^{\gamma} s^{k-1}ds\\
&=&\int_0^1 \prod_{i=1}^n\left(\sum_{\delta_i =0}^{\gamma_i} \frac{\gamma_i !}{\delta_i ! (\gamma_i - \delta_i)!}y_i^{(\alpha_i - \delta_i)}\Delta x^{\delta_i} s^{\delta_i}\right)s^{k-1}ds\nonumber\\
&=&\int_0^1\sum_{\delta=0}^{\gamma}\left(  \prod_{i=1}^n\frac{\gamma_i !}{\delta_i ! (\alpha_i - \delta_i)!}y_i^{(\gamma_i - \delta_i)}\Delta x^{\delta_i}\right) s^{\delta} s^{k-1} ds \nonumber\\
&=&\sum_{\delta=0}^{\gamma}\frac{1}{\delta_{\Sigma} +k}\left[  \prod_{i=1}^n\left( \frac{\gamma_i !}{\delta_i ! (\gamma_i - \delta_i)!}\right)  y^{(\gamma - \delta)}\right] \Delta x^{\delta}\nonumber\\
&=:&\sum_{\delta=0}^{\gamma}p_{k\delta}^{y\gamma} \Delta x^{\delta}\nonumber
\end{eqnarray}
where $\displaystyle \delta_{\Sigma} := \sum_{i=1}^n \delta_i$ and $s^{\delta}=\Pi_{i=1}^{n}s^{\delta_i}=s^{\delta_{\Sigma}}$. Hence
\begin{equation}
\begin{array}{ll}
 c^j_{k} (x,y) =\\
\\
 \sum_{\gamma}{\big \{}\sum_i\sum_{\beta + \alpha = \gamma} (\beta_i +1)(\alpha_i +1 )c^j_{r(\beta +1_i)}c^j_{(k-1-r)(\alpha +1_i)}   \\
\\
 +\sum_i (\gamma_i+2)(\gamma_i+1)c_{k(\gamma+2_i)}+ 
\sum_{\beta + \alpha = \gamma}(\sum\frac{1}{\beta !}b^{j}_{lm,\beta}(y)\times\\
\\
(\alpha_i+1)c_{(k-1)(\alpha+1_i)}{\big \}}\left( \sum_{\delta=0}^{\gamma}p_{k\delta}^{y\gamma} \Delta x^{\delta}\right).  
\end{array}
\end{equation}

\section{Proof of convergence of the formal solution (\ref{pjo}) for some time $0\leq t\leq T_0$}

In this section we shall show that the representation (\ref{pjo}) of the solution of (\ref{parasyst2}) holds for some time $0\leq t\leq T_0$ for some $T_0>0$. In the next section then we will show how the time transformation $t\rightarrow \tau(t)$ is used to get a representation of the fundamental solution for any domain with any finite time horizon. We shall prove that the representation (\ref{pjo}) holds for the equation
\begin{equation}\label{parasystbeta}
\frac{\partial u_i}{\partial \tau}=\beta\sum_{j=1}^n \frac{\partial^2 u_i}{\partial x_j^2} 
+\beta\sum_{j,k=1}^n b^i_{jk}\frac{\partial u_j}{\partial x_k}
\end{equation}
with some $\beta$ such that $t=\tau\beta$.
Essentially this step means that we have proved the validity of the representation for smaller time. Considering the solution of (\ref{parasyst2}) at time $t$ is equivalent to considering the solution of (\ref{parasystbeta}) at time $\tau$ with $t=\tau\beta$. If we want coefficients to be small then we choose $\beta$ small. Hence if 
\begin{equation}\label{pjobeta}
p^{\beta}_j(\tau,x,y)=\frac{1}{\sqrt{4\pi \tau}^n}\exp\left(\frac{\sum_{i=1}^n \Delta x_i^2}{4\tau}+\sum_{k=0}^{\infty} c^j_{k,\beta}(x,y)\tau^k \right)  
\end{equation}
is a representation of the fundamental solution of (\ref{parasystbeta}) on some domain $\Omega\times (0,T]$ for some $\tau <1$ (which may be shown by proving that for each $x,y\in \Omega$ $c^j_{k,\beta}( x,y)$ converges to zero), then this is a proof that the representation (\ref{pjo}) of the fundamental solution of (\ref{parasyst2}) converges for $t=\beta\tau$.

\subsection{Proof of convergence on bounded domains }
Since $\Omega\subset {\mathbb R}^n$ is bounded, there is a ball $B_R(0)$ around $0$ with radius $R$ such that $\Omega\subset B_R(0)$. Recall that
\begin{equation}\label{cj0}
 c^j_0(x,y)=\sum_m (y_m-x_m)\int_0^1 \sum_lb^j_{l,m}(y+s(x-y))ds,
\end{equation}
and for all $k\geq 1$ we have
\begin{equation}
c^j_k(x,y)=\int_0^1 R_{k-1}(y+s(x-y),y)s^{k-1}ds
\end{equation}
with 
\begin{equation}\label{tk}
\begin{array}{ll}
R^j_{k-1}(t,x,y)=&\Delta c^j_{k-1}+\sum_{l=1}^n\sum_{r=0}^{k-1}\left( \frac{\partial}{\partial x_l}c^j_r\frac{\partial}{\partial x_l}c^j_{k-1-r}\right)\\
\\
&+\sum_{lm} b^j_{lm}(x)\frac{\partial}{\partial x_m}c^l_{k-1}.
\end{array}
\end{equation}
If the modulus of the coefficients $b^j_{l,m}$ are bounded by the generic $C >0$ for all $j,l,m$, then we have
\begin{equation}\label{cjineq}
 |c^j_0(x,y)|\leq  n^2RC.
\end{equation} 
Next the time transformation
\begin{equation}
t=\beta \tau
\end{equation}
transforms the equation
\begin{equation}\label{parasyst3a}
\frac{\partial u_i}{\partial t}=\sum_{j=1}^n \frac{\partial^2 u_i}{\partial x_j^2} 
+\sum_{jk} b^i_{jk}\frac{\partial u_j}{\partial x_k}
\end{equation}
into the equation
\begin{equation}\label{parasyst3a} 
\frac{\partial v_i}{\partial \tau}=\beta\sum_{j=1}^n \frac{\partial^2 v_i}{\partial x_j^2} 
+\beta\sum_{jk} b^i_{jk}\frac{\partial v_j}{\partial x_k},
\end{equation}
where $u_i(t,x)=v_i(\tau,x)$, where $\frac{\partial u_i}{\partial t}=\frac{\partial v_i}{\partial \tau}\frac{\partial \tau}{\partial t}=\frac{\partial v_i}{\partial \tau}\frac{1}{\beta}$.
The analogous representation of the solution is of the form
\begin{equation}\label{betapj}
p^{\beta}_j(t,x,y)=\frac{1}{\sqrt{4\pi \tau}^n}\exp\left(\frac{\sum_{i=1}^n \Delta x_i^2}{4\beta\tau}\right)\left(1+\sum_k c^{j}_{k,\beta}(x,y)\tau^k \right),  
\end{equation}
for $j=1,\cdots ,n$. Plugging (\ref{betapj}) into (\ref{parasyst3a}) and collecting the terms with $\tau^{-2}$, $\tau^{-1}$ etc. we get (we feel free to write $t$ instead of $\beta \tau$ if convenient)
\begin{equation}\label{tauminus2}
\tau^{-2}:~~\frac{\sum_i \Delta x_i^2}{4\beta \tau^2}=\beta\sum_{l}\frac{ \Delta x_l^2}{4\beta^2 \tau^2}
\end{equation}

\begin{equation}\label{tauminus1}
\tau^{-1}:~~-\frac{n}{2t}=-\beta\sum_{l}\frac{1}{2\beta t}-\frac{\beta}{2\beta \tau}\left( \sum_l \Delta x_l\frac{\partial c_{0,\beta}^j}{\partial x_l}-\sum_{lm} b^j_{lm}(x)\Delta x_m\right) , 
\end{equation}
and for all $k-1\geq 0$
\begin{equation}\label{tk}
\begin{array}{ll}
\tau^{k-1}:~~kc^j_{k,\beta}+\beta\sum_l \Delta x_l\frac{\partial c_{k,\beta}^j}{\partial x_l}=&\beta\Delta c^j_{k-1,\beta}+\beta\sum_{l=1}^n\sum_{r=0}^{k-1}\left( \frac{\partial}{\partial x_l}c^j_{r,\beta}\frac{\partial}{\partial x_l}c^j_{k-1-r,\beta}\right)\\
\\
&+\beta\sum_{lm} b^j_{lm}(x)\frac{\partial}{\partial x_m}c^l_{k-1,\beta}\equiv \beta R^j_{k-1}(x,y).
\end{array}
\end{equation} 
We divide equation (\ref{tk}) by $\beta$ and get the
solutions (the solution for $c^{j}_{0,\beta}$ equals exactly that for $c^{j}_0$ in (\ref{cj0}))
\begin{equation}\label{jbk}
c^{j}_{k,\beta}(x,y)=\int_0^1 R^j_{k-1}(y+s(x-y),y)s^{\frac{k}{\beta}-1}ds.
\end{equation}
Next we prove 
\begin{thm}
There exists $\beta >0$ such that 
\begin{equation}
\mbox{for each}~~x,y\in \Omega, 1\leq j\leq n~~c^j_{k,\beta}(x,y)\downarrow 0 \mbox{ as } k\uparrow \infty.
\end{equation}
\end{thm}
\begin{proof}
First we remark that
\begin{equation}
\frac{\partial^{|\alpha|}}{\partial x^{\alpha}}c_0(x,y)\leq c^{|\alpha|}c^{up}_0,
\end{equation}
for some generic $C$, where
\begin{equation}
c^{up}_0:=\sup_{x,y\in \Omega}c_0(x,y).
\end{equation}
Indeed, if we define
\begin{equation}
b^{up}:=\sup_{x\in{\Omega},1\leq l,m\leq n}b^j_{lm}(x),
\end{equation}
then writing the multivariate derivative of order $\alpha$ with $\alpha=(\alpha_1,\cdots ,\alpha_n)$, and $|\alpha|:=\sum_{i=1}^n\alpha_i$ we get $|\alpha|n$ terms bounded by $C^{|\alpha|-1}b^{up}$ and $n^2$ terms bounded by $RC^{|\alpha|}b^{up}$.

Next, a majorant of $c^{j}_{k,\beta}(x,y)$ is obtained as follows: we consider three types of operators $O^{1,n}_k, O^{2,n}_k, O^{3,n}_k$ with positive integers $k$, and acting on a single function $f:\Omega\times \Omega\rightarrow {\mathbb R}$ or on a families of functions $(f_l)_{1\leq l\leq k}:\Omega\times \Omega\rightarrow {\mathbb R}$, namely
\begin{equation}
\begin{array}{ll}
O^{1,n}_k\left[f\right](x,y):=\frac{\beta}{k}\Delta f(x,y)\\
\\
O^{2,n}_k\left[f_k,\cdots,f_1\right](x,y):=\frac{\beta}{k}\sum_{l=1}^n\sum_{r=0}^k \frac{\partial f_r}{\partial x_l}\frac{\partial f_{k-r}}{\partial x_l}\\
\\
O^{3,n}_k\left[f\right](x,y):=\frac{\beta}{k}\sum_{lm} b^j_{lm}(x)\frac{\partial}{\partial x_m}f(x,y).
\end{array}
\end{equation}
Let
\begin{equation}
M_k:=\left\lbrace (\alpha_k,\cdots ,\alpha_1)|\alpha_j\in \left\lbrace 1,2,3\right\rbrace \right\rbrace 
\end{equation}
For 
\begin{equation}
c^{up}_{k,\beta}:=\sup_{x,y\in\Omega, j\in\left\lbrace 1,\cdots,n\right\rbrace}c^j_{k,\beta}(x,y) 
\end{equation}
we have
\begin{equation}\label{boundest}
c_{k,\beta}^{up}\leq \sum_{\alpha\in M_k}O_k^{\alpha ,n}c_0(x,y)=\sum_{\alpha\in M_k}O_k^{\alpha ,n}c_0(x,y),
\end{equation}
where
\begin{equation}
O^{\alpha ,n}_k\left[f\right](x,y):=O^{\alpha_k,n}_kO^{\alpha_{k-1},n}_{k-1}\circ\cdots\circ O^{\alpha_1,n}_1\left[f\right](x,y). 
\end{equation}
First let ${\bf 1}_k$ (resp. ${\bf 2}_k, {\bf 3}_k$) the multiindex $\alpha \in M_k$ such that for each $1\leq m\leq k$ $\alpha_m=1$ (resp. $\alpha_m=2, \alpha_m=3$). Hence
\begin{equation}\label{abb}
\begin{array}{ll}
O^{{\bf 1}}_k\left[f\right] (x,y)=\left( O^1\right)^k\left[f\right] (x,y)=\Delta^k\left[f\right] (x,y) 
\end{array}
\end{equation}
etc.. Then for all $x,y\in \Omega$
\begin{equation}
|O^{{\bf 1,n}}_k\left[c_0\right] (x,y)|\leq \frac{\beta^kn^{k}C^{2k}c_0^{up}}{k!},
\end{equation}
and for $b$ bounded by a generic $C$ on the domain we also have
\begin{equation}
|O^{{\bf 3,n}}_k\left[c_0\right] (x,y)|\leq \frac{\beta^k n^{2k}C^{2k}c_0^{up}}{k!}.
\end{equation}
The operators of quadratic type applied to $c_0(x,y)$ $O^{2,n}_{k}c_0$ decrease to zero as $k\uparrow \infty$ if $\beta$ is small. We estimate (a rough estimate is sufficient here)
that surely we have
\begin{equation}\label{large}
|O^{{\bf 2},n}_k\left[f\right] (x,y)|\leq \frac{\beta^k k^2 2^{k-1}n^kC^{k}k!(c_0^{up})^{k+1}}{k!}.
\end{equation}
For large $k$ this is essentially the largest term of all the $3^k$ contributions in the sum (\ref{boundest}) for large $k$ ($k$ fixed). We mean the following: if we choose 
\begin{equation}\label{beta}
\beta <\frac{1}{3\cdot 4n^2C^2 (c_0^{up})^2},
\end{equation}
then surely we have for $k\geq k_0$ (some $k_0>0$
\begin{equation}\label{essest}
|3^kO^{{\bf 2}}_k\left[c_0\right] (x,y)|\leq \frac{3^k\beta^k k^2 2^{k-1}n^kC^{k}k!(c_0^{up})^{k+1}}{k!}\downarrow 0.
\end{equation}
as $k\uparrow \infty$, and this is also the estimate which holds for $c_k$ for large $k$. Here we choose $\beta$ such that in a summand in $O_k^{\alpha ,n}c_0(x,y)$ in (\ref{boundest}) each occurrence of an operator of form $O^{3,n}_k$ can be replaced by an operator of form  $O^{2,n}_k$ in order to get a majorant estimation. So in the sum in (\ref{boundest}) it suffices to concentrate on the summands consisting of concatenations of operators of form $O^{2,n}_k$ and $O^{1,n}_k$. For natural numbers $l$ let us define an increasing sequence of numbers $k_1<k_2<\cdots <k_l<k_{l+1}\cdots$, and operators
\begin{equation}
\begin{array}{ll}
O^{1,n}_{k_{l+1}k_l}:=O^{1,n}_{k_{l+1}}\circ\cdots \circ O^{1,n}_{k_l}\\
\\
O^{2,n}_{k_{l+1}k_l}:=O^{2,n}_{k_{l+1}}\circ\cdots \circ O^{2,n}_{k_l}
\end{array}
\end{equation}
Then in the summands o (\ref{boundest}) we have to consider the asymptotic behavior of values of family of operators of form
\begin{equation}\label{fam1}
O^{2,n}_{k_{l+1}k_l}\circ O^{1,n}_{k_{l}k_{l-1}}\circ \cdots \circ O^{2,n}_{k_{3}k_2}\circ O^{1,n}_{k_{2}k_1}
\end{equation}
or of form
\begin{equation}\label{fam2}
O^{2,n}_{k_{l+1}k_l}\circ O^{1,n}_{k_{l}k_{l-1}}\circ \cdots \circ O^{1,n}_{k_{3}k_2}\circ O^{2,n}_{k_{2}k_1}
\end{equation}
applied to $c_0(x,y)$ as $k\uparrow \infty$. If there is only a finite occurrence of operators of form $O^{1,n}_k$ in such a family ((\ref{fam1}) or (\ref{fam2})), then the asymptotic behavior is clearly the same as for 
$O^{{\bf 2}}_k c_0(x,y)$. 
If on the other hand there are infinite occurrences of operators of form $O^{1,n}_k$ 
in ((\ref{fam1}) of (\ref{fam2})), then for large $k$ $O^{{\bf 2}}_k c_0(x,y)$ becomes a 
majorant of such a summand. Hence, the estimate (\ref{essest}) is a majorant for large $k$ and proves the convergence of the series in (\ref{boundest}). 
\end{proof}

\subsection{Remark on unbounded domains}

It is not possible to extend the proof in the preceding section immediately to unbounded domains $\Omega \subseteq {\mathbb R}^n$. However, a similar proof with a different but equivalent ansatz 
\begin{equation}\label{pdordj}
p^{d}_j(t,x,y)=\frac{1}{\sqrt{4\pi t}^n}\exp\left(-\frac{\sum_{i=1}^n \Delta x_i^2}{4t}\right)\left( 1 +\sum_{k=0}^{\infty} d^j_{k}(x,y)t^k \right),  
\end{equation}
leads to such an extension. The recursion equation for $d_0$ and $c_0$ are equivalent, but the recursion equations for the $d_k, k\geq 1$ are more involved. However, it can be shown that given $t,y$ the supremum in $x$ of each
\begin{equation}\label{wc0}
\frac{1}{\sqrt{4\pi t}^n}\exp\left(\frac{-\Delta x^2}{4t} \right)d^j_k(x,y)
\end{equation}
is in some ball which can be chosen a priori. However, this is beyond the scope of the present paper, and we shall consider a similar situation in \cite{KG}.

\section{Extension to the time-inhomogeneous case (solution and global convergence)}
In a second step we use a certain nonlinear time transformation in order to lift the result to convergence for any finite time $t$. This requires the extension of the analysis to the case with time-dependent coefficients. Note that in the extension of the recursion of the $c_k$ to the time-inhomogeneous case only first order time derivatives occur. This is the reason for the weaker constraints for (\ref{coeft}).
We start this Section with the computation of the recursive coefficients $c_k$ in the case of time- and space-dependent drift coefficients $b^j_{kl}$. Then we shall complete the proof for convergence on bounded domains for any finite time in the time-homogenous case, and finally in the time-inhomogeneous case in the following subsections. 

\subsection{Formal computation of recursive coefficients in the time-inhomogeneous case}
We consider parabolic equations with time-dependent coefficients of the form
\begin{equation}
\frac{\partial u^i}{\partial t}+\Delta u^i +\sum_{jk} b^i_{jk}(t,x)\frac{\partial u^j}{\partial x_k}=0
\end{equation}
We consider the ansatz
\begin{equation}
p_j(t,x,0,y)=\frac{1}{\sqrt{4\pi t}^n}\exp\left(-\frac{\Delta x^2}{4t}+ \sum_{k=0}^{\infty}c^j_k(t,x,y)t^k\right).
\end{equation}
Compared to the time-homogenous case the time derivative contains an additional term. We have
\begin{equation}
\begin{array}{ll}
\frac{\partial p_j}{\partial t}(t,x,y)={\Bigg (}-\frac{n}{2t}+\frac{\sum_i \Delta x_i^2}{4t^2}+\sum_{k=0}^{\infty}\frac{\partial c_k}{\partial t}(t,x,y)t^k\\
\\
\hspace{5.5cm}+\sum_k kc^j_k(t,x,y)t^{k-1}{\Bigg )}p_j (t,x,0,y)
\end{array}
\end{equation}
The spatial derivatives are essentially the same as in the time-homogenous case. We compute
\begin{equation}
\frac{\partial p_j}{\partial x_l}(t,x,y)=\left(\frac{-\Delta x_l}{2t}+\sum_k \frac{\partial}{\partial x_l}c^j_k(t,x,y)t^k \right) p_j(t,x,0,y),
\end{equation}
and
\begin{equation}
\begin{array}{ll}
\frac{\partial^2 p_j}{\partial x_l^2}(t,x,y)=&{\Bigg (}-\frac{1}{2t}+\sum_k \frac{\partial^2}{\partial x_l^2}c^j_k(t,x,y)t^k \\
\\
&+\left(-\frac{\Delta x_l}{2t}+\sum_k \frac{\partial}{\partial x_l}c^j_k(t,x,y)t^k \right)^2{\Bigg )} p_j(t,x,0,y).
\end{array}
\end{equation}
Plugging into \ref{parasyst1} and ordering with respect to the terms $t^{-2},t^{-1}$ etc. we get the following
recursive relations for the $c^j_k$, where $1\leq j\leq n$. First, the highest order terms are the same as before:
\begin{equation}\label{tminus2}
t^{-2}:~~\frac{\sum_i \Delta x_i^2}{4t^2}=\sum_{l}\frac{ \Delta x_l^2}{4t^2}
\end{equation}
The terms of order $t^{-1}$ are essentially as before (we just have to add the $t$-argument in the coefficient functions $b^i_{jk}$):
\begin{equation}\label{tminus1}
t^{-1}:~~-\frac{n}{2t}=-\sum_{l}\frac{1}{2t}-\frac{1}{2t}\left( \sum_l \Delta x_l\frac{\partial c_0^j}{\partial x_l}-\sum_{lm} b^j_{lm}(t,x)\Delta x_m\right). 
\end{equation}
For $k-1\geq 0$ we get an additional $t$-derivative on the right side:
\begin{equation}\label{timetk}
\begin{array}{ll}
t^{k-1}:~~kc^j_k+\sum_l \Delta x_l\frac{\partial c_k^j}{\partial x_l}=\frac{\partial c^j_{k-1}}{\partial t}+\Delta c^j_{k-1}+\sum_{l=1}^n\sum_{r=0}^{k-1}\left( \frac{\partial}{\partial x_l}c^j_r\frac{\partial}{\partial x_l}c^j_{k-1-r}\right)\\
\\
+\sum_{lm} b^j_{lm}(t,x)\frac{\partial}{\partial x_m}c^l_{k-1}\equiv R^j_{k-1}(x,y)
\end{array}
\end{equation}
Hence,
\begin{equation}
\sum_l \Delta x_l\frac{\partial c_0^j}{\partial x_l}=-\sum_{l,m} b^j_{lm}(t,x)\Delta x_m, 
\end{equation}
which has the solution
\begin{equation}
 c^j_0(x,y)=\sum_m (y_m-x_m)\int_0^1 \sum_lb^j_{l,m}(t,y+s(x-y))ds,
\end{equation}
and for all $k\geq 1$ we have
\begin{equation}
c^j_k(x,y)=\int_0^1 R_{k-1}(t,y+s(x-y),y)s^kds
\end{equation}
with $R_{k-1}$ as in equation \eqref{timetk}.
The explicit calculation of the solution is know completely analogous, so it suffices to write down the results.
We write
\begin{equation}
b^j_{lm}(t,x)=\sum_{\gamma} \frac{1}{\gamma !}b^{j}_{lm,\gamma}(t,y)(\Delta x)^{\gamma}
\end{equation}  
along with some multiindex $\gamma$. Then
\begin{equation}
\begin{array}{ll}
 c^j_0(t,x,y)&=-\sum_{l,m} \sum_{\gamma}b^{j}_{lm\gamma}(y)\Delta x^{\gamma+1_i}\frac{1}{1+|\gamma|}\\
\\
&\equiv \sum_{\gamma} c^j_{0\gamma}(t,y)\Delta x^{\gamma}
\end{array}
\end{equation}
Given that $c^ j_{k-1}$ equals its Taylor series for every $y\in {\mathbb R}^n$, i.e.
\begin{equation}
c^j_{k-1}(t,x)=\sum_{\gamma} c^j_{(k-1)\gamma}(t,y)\Delta x^{\gamma}=\sum_{\gamma,l}c^j_{(k-1)\gamma l}(y)\Delta x^{\gamma}t^l,
\end{equation}
we have
\begin{equation}
\begin{array}{ll}
 c^j_{k} (t,x,y) =\sum_{\gamma,l}lc^j_{(k-1)\gamma l}(y)\Delta x^{\gamma}t^l\\
\\
+\sum_{\gamma}{\big \{}\sum_i\sum_{\beta + \alpha = \gamma} (\beta_i +1)(\alpha_i +1 )c^j_{r(\beta +1_i)}(t,y)c^j_{(k-1-r)(\alpha +1_i)}(t,y)   \\
\\
 +\sum_i (\gamma_i+2)(\gamma_i+1)c_{k(\gamma+2_i)}+ 
\sum_{\beta + \alpha = \gamma}(\sum\frac{1}{\beta !}b^{j}_{lm,\beta}(t,y)\times\\
\\
(\alpha_i+1)c_{(k-1)(\alpha+1_i)}{\big \}}\left( \sum_{\delta=0}^{\gamma}p_{k\delta}^{y\gamma} \Delta x^{\delta}\right),  
\end{array}
\end{equation}
where the $p_{k\delta}^{y\gamma}$ are defined exactly as before. 

\subsection{Completion of convergence proof for finite time in the general case on bounded  domains}

We apply a nonlinear time transformation. Consider for $\beta >0$ the transformation $\tau(t):[0,\infty)\rightarrow [0,1)$
\begin{equation}
\tau =(1-e^{-\frac{t}{\beta}}),~~\mbox{or}~~t=t(\tau)=-\beta\ln(1-\tau)
\end{equation}
with
\begin{equation}
\frac{\partial t}{\partial \tau}=\frac{\beta}{1-\tau}.
\end{equation}
The transformation of the equation 
\begin{equation}\label{tb}
\frac{\partial u_i}{\partial t}=\Delta u_i +\sum_{jk}b^i_{jk}(t,x)\frac{\partial u_j}{\partial x_k}
\end{equation}
then is
\begin{equation}\label{taub}
\frac{\partial v_i}{\partial \tau}=\frac{\beta}{1-\tau}\Delta v_i +\frac{\beta}{1-\tau}\sum_{jk}b^i_{jk}(t(\tau),x)\frac{\partial v_j}{\partial x_k}.
\end{equation}
Let us call the associated coefficients of the global expansion of the fundamental solution by $c^{j}_{k,\beta,\tau}$. If we can show that for each given $x,y$ (in $\Omega$ and then in ${\mathbb R}^n$ in general) $c^{j}_{k,\beta,\tau}(\tau,x,y)$ converges to zero as $k\uparrow \infty$, then we have convergence for $\tau<1$ which implies convergence of the analytic representation for the original equation for
$t\in (0,\infty)$.  First we derive the recursive relations for (\ref{taub}). Since
\begin{equation}
\phi^{\beta,\tau}_i(\tau,x,y):=\frac{1}{\sqrt{4\pi\left(-\beta\ln(1-\tau)\right) }^n}\exp\left(-\frac{\Delta x^2}{4\left(-\beta \ln(1-\tau) \right) }\right) 
\end{equation}
is the fundamental solution of the equation
\begin{equation}\label{ttau}
\frac{\partial u}{\partial \tau}=\frac{\beta}{1-\tau}\Delta u,
\end{equation}
we consider the ansatz
\begin{equation}
p^{\beta,\tau}_i(\tau,x,y)=\phi^{\beta,\tau}_i(\tau,x,y)\exp\left(\sum_{k=0}^{\infty}c^i_{k,\beta,\tau}(\tau,x,y)\tau^k\right).
\end{equation}
Using $t=-\beta\ln(1-\tau)$ we have
\begin{equation}
\begin{array}{ll}
\frac{\partial p_i}{\partial \tau}(t,x,y)={\Bigg (}-\frac{n}{2t}\frac{\partial t}{\partial \tau}+\frac{\sum_i \Delta x_i^2}{4t^2}\frac{\partial t}{\partial \tau}+\sum_{k=0}^{\infty}\frac{\partial }{\partial \tau}c^i_{k,\beta, \tau}(\tau,x,y)\tau^k\\
\\
\hspace{5.5cm}+\sum_k kc^i_{k,\beta,\tau}(\tau,x,y)\tau^{k-1}{\Bigg )}p^{\beta,\tau}_i (\tau,x,y),
\end{array}
\end{equation}
\begin{equation}
\frac{\partial p_i}{\partial x_l}(\tau,x,y)=\left(\frac{-\Delta x_l}{2t}+\sum_k \frac{\partial}{\partial x_l}c^j_{k,\beta,\tau})(\tau,x,y)\tau^k \right) p^{\beta,\tau}_i(\tau,x,y),
\end{equation}
and
\begin{equation}
\begin{array}{ll}
\frac{\partial^2 p_i}{\partial x_l^2}(\tau,x,y)=&{\Bigg (}-\frac{1}{2t}+\sum_k \frac{\partial^2}{\partial x_l^2}c^j_{k,\beta,\tau}(\tau,x,y)\tau^k \\
\\
&+\left(-\frac{\Delta x_l}{2t}+\sum_k \frac{\partial}{\partial x_l}c^j_{k,\beta,\tau}(\tau,x,y)\tau^k \right)^2{\Bigg )} p^{\beta,\tau}_i(\tau,x,y).
\end{array}
\end{equation}
Plugging into (\ref{taub}) and ordering with respect to the terms $\tau^{-2},\tau^{-1}$ etc. leads to
\begin{equation}\label{taubetaminus2}
\tau^{-2}:~~\frac{\sum_i \Delta x_i^2}{4t^2}\frac{\partial t}{\partial \tau}=\frac{\beta}{1-\tau}\sum_{l}\frac{ \Delta x_l^2}{4t^2},
\end{equation}
which is satisfied because the second order diffusion term in (\ref{taub}) is $\frac{\beta}{1-\tau}$. 
For the terms of order $\tau^{-1}$ we get: 
\begin{equation}\label{tbetaminus1}
\begin{array}{ll}
\tau^{-1}:~~-\frac{n}{2t}\frac{\partial t}{\partial \tau}=-\frac{\beta}{1-\tau}\sum_{l}\frac{1}{2t}\\
\\
-\frac{1}{2t}\left( \frac{\beta}{1-\tau}\sum_l \Delta x_l\frac{\partial c_0^j}{\partial x_l}-\frac{\beta}{1-\tau}\sum_{lm} b^j_{lm}(t,x)\Delta x_m\right). 
\end{array}
\end{equation}
For $k-1\geq 0$ we get an additional $\tau$-derivative on the right side:
\begin{equation}\label{timebetatk}
\begin{array}{ll}
\tau^{k-1}:~~kc^i_k+\frac{\beta}{1-\tau}\sum_l \Delta x_l\frac{\partial c_k^i}{\partial x_l}=\frac{\partial c^i_{k-1}}{\partial \tau}+\frac{\beta}{1-\tau}\Delta c^i_{k-1}\\
\\
+\frac{\beta}{1-\tau}\sum_{l=1}^n\sum_{r=0}^{k-1}\left( \frac{\partial}{\partial x_l}c^i_r\frac{\partial}{\partial x_l}c^i_{k-1-r}\right)\\
\\
+\frac{\beta}{1-\tau}\sum_{lm} b^j_{lm}(t,x)\frac{\partial}{\partial x_m}c^l_{k-1}\equiv \frac{\beta}{1-\tau}R^i_{k-1}(\tau,x,y)
\end{array}
\end{equation}
We have
\begin{equation}
 c^i_{0,\beta,\tau}(\tau,x,y)=\sum_m (y_m-x_m)\int_0^1 \sum_lb^i_{l,m}(t(\tau),y+s(x-y))ds,
\end{equation}
and for all $k\geq 1$ we have 
\begin{equation}
c^i_{k,\beta,\tau}(x,y)=\int_0^1 R^{i,\tau}_{k-1}(t(\tau),y+s(x-y),y)s^{\frac{(1-\tau)k}{\beta}-1}ds,
\end{equation}
where
\begin{equation}
R^{i,\tau}_{k-1}(\tau,x,y)=\frac{1-\tau}{\beta}\frac{\partial }{\partial \tau}c^i_{k-1,\beta,\tau}+R^{i,\tau}_{k-1,\beta,\tau}(\tau,x,y)
\end{equation}
with $R^i_{k-1}$ is as in equation \eqref{timebetatk}. From this and the preceding sections it is clear how we get the power series representation (\ref{powerck})in theorem 2.1. above.
We see from this representation that the proof for small $t$ can be used, only that the substitution
\begin{equation}
\beta \rightarrow \frac{\beta}{1-\tau}
\end{equation}
has to be made. Since there are only first order time derivatives in the recursion (cf. (\ref{timetk}) and (\ref{timebetatk}), the proof of section 5.1. can be extended trivially. Hence, global convergence (for any positive $t$ of our analytic expansion follows from the following 
\begin{prop}
For each constant $c>0$ the range of the function 
\begin{equation}
(\beta ,\tau)\rightarrow t=-\beta\ln(1-\tau)
\end{equation}
is unbounded on the domain
\begin{equation}
\left\lbrace (\beta,\tau)| \frac{\beta}{1-\tau}=c\right\rbrace .
\end{equation}
\end{prop}
\begin{proof} $c=\frac{\beta}{1-\tau}=\frac{\epsilon\beta}{\epsilon(1-\tau)}\rightarrow -\epsilon\beta\ln(\epsilon(1-\tau))\uparrow \infty$ as $\epsilon \downarrow 0$.\end{proof}
This means that it suffices to prove that the recursion (\ref{tbetaminus1}, (\ref{timebetatk}) converges to zero for some $\frac{\beta}{1-\tau}$ (which may be as small as we need).

\section{Representations of initial boundary value problems of first and second type}

The explicit fundamental solution leads to representations of solutions for  initial-boundary problems of parabolic systems and parabolic equations. We consider two examples. 
\subsection{Representation of the solution for initial-boundary problems for parabolic systems of first type}
For the Cauchy problem (\ref{cp}) we have the following representation of the solution $u$:
\begin{equation}
\begin{array}{ll}
u(t,x)=\int_{{\mathbb R}^n}\phi(y){\mathbf p}(t,x;0,y){\mathbf \phi}(y)dy\\
\\
+\int_0^t\int_{{\mathbb R}^n}f(s,y){\mathbf p}(t,x;s,y){\mathbf \phi}(y)dyds
\end{array}
\end{equation}
\begin{rem}
Strictly speaking, the solution for ${\bf p}$ presented here is on bounded domains $\Omega \subset {\mathbb R}^n$ (which is for large $\Omega$ a sufficient approximation for numerical purposes, but not exact). However, such exact representations on unbounded domains can be found using the recursion indicated in Section 4.2. 
\end{rem}

\subsection{Representation of the solution for initial-boundary problems for parabolic equations of second type}

In the case of the scalar problem \ref{ibst} for the solution $u$ the ansatz for $\gamma$ with
\begin{equation}
\begin{array}{ll}
u(t,x)=&\int_{\Omega}\phi(y)p(t,x,0,y)dy-\int_0^t\int_{\Omega}f(s,y)p(t,x,0,y)dyds\\
\\
&+\int_0^t\int_B p(t,x;s,y)\gamma (s,y)ds dy
\end{array}
\end{equation}
leads to the integral equation
\begin{equation}
\begin{array}{ll}
\frac{1}{2}\gamma(t,x)=&\int_0^t\int_B\left\lbrace \frac{\partial p}{\partial \nu}(t,x;s,y)+\alpha(t,x)p(t,x;s,y)\right\rbrace \gamma(t,x)dB_xds\\
\\
& +h(t,x)
\end{array}
\end{equation}
where 
\begin{equation}
\begin{array}{ll}
h(t,x)=\int_{\Omega}\frac{\partial p}{\partial \nu(t,x)}(t,x;s,y)\phi(y)dy\\
\\
+\int_0^t\int_{\Omega}\frac{\partial p}{\partial \nu(t,x)}(t,x;s,y)f(s,y)dyds\\
\\
+\alpha(t,x)\int_{\Omega}p(t,x;s,y)\phi(y)dy\\
\\
-\alpha(t,x)\int_0^t\int_{\Omega}p(t,x;s,y)f(s,y)dyds\\
\\
-\psi(t,x)
\end{array}
\end{equation}
Hence with our explicit solution for $p$ we reduce the initial-boundary value problem of second type to the solution of a linear integral equation.

\section{Generalizations, applications, and further comments}

The preceding results can be extended to more general diffusions. 
We have
\begin{thm}
Consider equation (\ref{parasyst2}) with
space-dependent diffusion coefficients $x\rightarrow a^i_{jk}(x)$ which satisfy
\begin{equation}\label{difcoeff}
|\partial_x^{\alpha}a^i_{jk}|\leq c^{|\alpha|}
\end{equation}
Assume that the conditions of theorem 1 are satisfied.
Then the fundamental solution has the representation
\begin{equation}\label{diffparasyst}
p^i(t,x,y)=\frac{1}{\sqrt{4\pi t}^n}\exp\left(-\frac{d_i^2(x,y)}{4t} \right)\exp\left(\sum_k c^i_{k,\beta}(\tau,x,y)\tau^k \right)  
\end{equation}
where for each $i$ $(x,y)\rightarrow d_i^2(x,y)$ are functionals which assign to each pair of points $x,y$ the
length of a geodesic with respect to the line element
\begin{equation}\label{line}
ds^2_i=\sum_{jk}g^i_{jk}dx_jdx_k,
\end{equation}
with $(g^i_{jk})$ the inverse of $(a^i_{jk})$, and the $c^i_k$ are smooth functions given by recursive relations similar to that in theorem 1 but involving $d^2_i$ and partial derivatives of $d^2_i$. 
\end{thm}
The proof is quite analogous except that additional existence and regularity results for the Riemmanian metric functional $d^2$ are needed. These are given in \cite{K2}. General analytical formulas are not available for the functional $d^2$ but in \cite{K2} that solutions can be approximated in arbitrarily strong Sobolev norms. This may be used to obtain approximations of (\ref{diffparasyst}) in arbitrarily strong Sobolev norms when combined with the results in \cite{K3}. Note, however, that an extension is far from obvious if the second order terms are coupled.

An immediate application of theorem 7.1. is a result of Varadhan which we state and prove here in the case of time-homogeneous coefficients and for scalar equations, where the highest order coefficient function in the global expansion may be denoted by $d^2$ without an index $i$.
\begin{cor}(time-homogeneous and scalar case)
Assume that for each $i$ we have $\lambda\xi^2\leq a^i_{jk}(x)\xi_i\xi_j \leq \Lambda \xi^2$ for $x\in \Omega\subseteq {\mathbb R}^n$ and some constants $0<\lambda <\Lambda$. Then for bounded H\"older-continuous coefficient functions $x\rightarrow a^i_{jk}(x)$, $x\rightarrow b^i_{jk}(x)$ 
\begin{equation}
\lim_{t\downarrow 0}-4 t \ln p(t,x,y) =d^2(x,y)
\end{equation}
where $d^2$ is the Riemannian metric induced by the line element (\ref{line}).
\end{cor}
\begin{proof}
The reason for the assumption of H\"older continuity is just for the existence of the fundamental solution, which may then be ensured by the parametrix method).
For the assumptions of theorem 7.1 this follows directly from the representation (\ref{diffparasyst}). Given $x,y$ one may define in a bounded domain $x,y\in \Omega$ containing the geodesic a series of coefficient functions $(a^{i,n}_{jk})_n$ and $(b^{i,n}_{jk})_n$ satisfying the assumptions of theorem 7.1. and such that $a^{i,n}_{jk}(x)\rightarrow a^i_{jk}(x)$ and $(b^{i,n}_{jk})_n\rightarrow b^i_{jk}$. Here we can assume that the corresponding geodesics connecting $x$ and $y$ are in $\Omega$ 
\end{proof}
There is a deep difference between the representations considered here with leading terms of the form
\begin{equation}
\frac{1}{\sqrt{4\pi t}^n}\exp\left(-\frac{d_i^2(x,y)}{4t} \right)
\end{equation}
and direct Taylor expansions of the solution. Indeed, in \cite{BKS} we saw that for the characteristic function (the Fourier transform of the fundamental solution with respect to the parameter $y$), where a direct Taylor approach seems natural, it seems that convergence results can be obtained only if coefficients are of linear spatial dependence.
We note that results of myself for scalar equations cited in \cite{KKS} cannot be directly generalized to the systems considered here. 
Our results  may be used to generalize the results in \cite{BS} and construct efficient computation schemes for related reaction diffusion equations. A second application may be the definition of generalized Brownian motions (cf. \cite{W}). This was attempted in \cite{R} in the context of elasticity and the Lam\'e equation, but not in a rigorous way. Note that Lam\'e equation has coupling of second order terms, so the generalized Brownian motions associated to (\ref{diffparasyst}) would not cover these examples from elasticity (because we have no second order coupling in (\ref{diffparasyst})). However, the functional analytic procedure to introduce processes as measures on path spaces using Riesz representation theorem leads to a new class of processes.
In the special case of higher dimensional scalar equations expansions of the type considered here in a probabilistic setting have been found to be very competitive (cf. \cite{KKS}). The results presented here are also a first step to get into deeper analysis of quasilinear parabolic systems, and numerical methods considered in (\cite{K3}) and (\cite{K2}) may be extended and used together with analytical results in (\cite{Kat}) and (\cite{L}) to obtain efficient and accurate schemes for quasilinear systems.

\section{Two examples of application: multidimensional Burgers system and Pauli equation}

The representation of solutions for parabolic systems considered in this paper has a wide range of applications.
One simple example is the multidimensional viscous Burgers equation with forcing. It is of the form
\begin{equation}
\frac{\partial  \mathbf{v}}{\partial t}-\nu \nabla^2 \mathbf{v}+ \left( \mathbf{v}\cdot\nabla\right)  \mathbf{u}=-\nabla F(t,x),
\end{equation}
where $x\in \Omega \subseteq {\mathbb R}^n$ is some domain, $F$ is some outer force, and $\mathbf{v}$ describes some velocity field. If the initial condition are of potential form, i.e.  $\mathbf{v}(0,x)=-\nabla \Phi(0,x)$ for some potential $\Phi$, then the velocity field remains a potential as time goes by and the dynamics id governed by 
\begin{equation}
\partial_t \Phi -\nu\nabla^2\Phi-\frac{1}{2}\nabla^2\Phi=F.
\end{equation}
Applying Cole-Hopf transformation $\psi=2\nu\ln \psi (t,x)$ we get an (imaginary-time) Schr\"odinger equation of form
\begin{equation}
\partial_t \psi -\nu\nabla^2\psi-\frac{1}{2\nu}F\psi.
\end{equation}
It is clear that our main results can be easily extended to this case with potential term.

Another example of interest is Pauli' s equation which describes the daynamics of on electron (with spin) in the presence of a (possibly space and time-dependent) magnetic field:
\begin{equation}
\begin{array}{ll}
ih\frac{\partial}{\partial t}\binom{\psi_1(t,x)}{\psi_2(t,x)}=&{\Bigg [}\left( \frac{-h^2}{2m_e}\mathbf{\Delta}+\frac{\mu_B}{h}\mathbf{L}\cdot \mathbf{B}+H_{\mbox{dia}} \right)  \mathbf{1}_2\\
\\
&+\mu_B\mathbf{\sigma}\cdot \mathbf{B}{\Bigg ]}\binom{\psi_1(t,x)}{\psi_2(t,x)},
\end{array}
\end{equation}
where $\mathbf{L}$ is the angular momentum $\mathbf{B}$ the magnetic field (possibly dependent on space and time), and $\mathbf{\sigma}$ is the vector of Pauli's spin matrices. Furthermore, as usual $m_e$ and $e$ denote the mass of the electron and the electric charge respectively, $h$ is the normalized Planck constant, and 
\begin{equation}
H_{\mbox{dia}}=\frac{e^2 \mathbf{B}^2}{8m_e}r^2\sin(\alpha)
\end{equation}
with $\alpha=\angle(x,\mathbf{B})$.
Again a slight extension of our proof of Theorem 2.1 leads to a global analytic expansion of our main result.

\end{document}